\documentclass[11pt]{article}
\usepackage{amsfonts}
\usepackage{amsmath}
\usepackage{amssymb}
\usepackage{graphicx}
\usepackage[dvipsnames,usenames]{color}
\usepackage[pdfstartview=FitH ]{hyperref}
\usepackage[hyperpageref]{backref}
\begin{document}

\baselineskip 18pt
\def\o{\over}
\def\e{\varepsilon}
\title{\Large\bf Mean\ Value\ from\ Representation\ of\ Rational\\
Number\ as\ Sum\ of\ Two\ Egyptian\ Fractions}
\author{Chaohua\ \  Jia}
\date{}
\maketitle {\small \noindent {\bf Abstract.} For given positive
integers $n$ and $a$, let $R(n;\,a)$ denote the number of positive
integer solutions $(x,\,y)$ of the Diophantine equation
$$
{a\o n}={1\o x}+{1\o y}.
$$
Write
$$
S(N;\,a)=\sum_{\substack{n\leq N\\ (n,\,a)=1}}R(n;\,a).
$$
Recently Jingjing Huang and R. C. Vaughan proved that for $4\leq N$
and $a\leq 2N$, there is an asymptotic formula
$$
S(N;\,a)={3\o \pi^2a}\prod_{p|a}{p-1\o p+1}\cdot N(\log^2N+c_1(a)
\log N+c_0(a))+\Delta(N;\,a).
$$
In this paper, we shall get a more explicit expression with better
error term for $c_0(a)$. }

\vskip.4in
\noindent{\bf 1. Introduction}

Representation of rational number as sum of unit fractions, or
Egyptian fractions, is an interesting topic in number theory. For
its history and related problems, one can see R. K. Guy's book[1].

Recently Jingjing Huang and R. C. Vaughan[2] studied the
representation of rational number as sum of two Egyptian fractions.
They established two mean value theorems, one of which is

\noindent{\bf Proposition 1.} For given positive integers $n$ and
$a$, let $R(n;\,a)$ denote the number of positive integer solutions
$(x,\,y)$ of the Diophantine equation
$$
{a\o n}={1\o x}+{1\o y}.
$$
Write
\begin{align*}
S(N;\,a)=\sum_{\substack{n\leq N\\ (n,\,a)=1}}R(n;\,a). \tag 1
\end{align*}
Then for $4\leq N$ and $a\leq 2N$, there is an asymptotic formula
$$
S(N;\,a)={3\o \pi^2a}\prod_{p|a}{p-1\o p+1}\cdot N(\log^2N+c_1(a)
\log N+c_0(a))+\Delta(N;\,a),
$$
where
\begin{align*}
c_1(a)=6\gamma-4{\zeta'(2)\o \zeta(2)}-2+\sum_{p|a}{6p+2\o p^2-1}
\cdot\log p \tag 2
\end{align*}
and
\begin{align*}
c_0(a)=-2\log^2a-4\log a\sum_{p|a}{\log p\o p-1}+O\Bigl({a\o
\varphi(a)}\cdot\log a\Bigr),\tag 3
\end{align*}
and
\begin{align*}
\Delta(N;\,a)\ll N^{1\o 2}\log^5N\cdot{a\o \varphi(a)}\prod_{p|a}
\Bigl(1-{1\o p^{1\o 2}}\Bigr)^{-1}. \tag 4
\end{align*}
Here $p$ denotes prime number, $\gamma$ is the Euler constant and
$\varphi(a)$ is the Euler totient function.

In this paper, we shall apply results in [3] and [4] to get a more
explicit expression with better error term for $c_0(a)$ in (3). We
shall prove

\noindent{\bf Theorem.} Let $S(N;\,a)$ be defined in (1). Then for
$4\leq N$ and $3\leq a\leq 2N$, we have
$$
S(N;\,a)={3\o \pi^2a}\prod_{p|a}{p-1\o p+1}\cdot N(\log^2N+c_1(a)
\log N+c_0(a))+\Delta(N;\,a),
$$
where $c_1(a)$ and $\Delta(N;\,a)$ are same as in (2) and (4), while
\begin{align*}
c_0(a)&=-2\log^2a-4\log a\sum_{p|a}{\log p\o p-1}-2\Bigl(\sum_{p|a}
{\log p\o p-1}\Bigr)^2\\
&\ +\Bigl(\sum_{p|a}{3p+1\o p^2-1}\cdot\log p\Bigr)^2-\sum_{p|a}
{3p^3+2p^2+3p\o (p^2-1)^2}\cdot\log^2 p
\end{align*}
\begin{align*}
&\ +\Bigl(6\gamma-4{\zeta'(2)\o \zeta(2)}-2\Bigr)\sum_{p|a}{3p+1\o
p^2-1}\cdot\log p \tag 5\\
&\ +{2\o \zeta^2(2)}(2\zeta'(2)+\zeta(2))^2-{4\o \zeta(2)}
(\zeta''(2)+\zeta'(2))\\
&\ -6\gamma\Bigl(2{\zeta'(2)\o \zeta(2)}+1\Bigr)+8\gamma^2-2\gamma_1
+2\zeta(2)+O\Bigl({1\o \varphi(a)}\Bigr).
\end{align*}
Here $p$ denotes prime number, $\varphi(a)$ is the Euler totient
function, $\gamma$ is the Euler constant and
\begin{align*}
\gamma_1=\lim_{H\rightarrow\infty}\Bigl(\sum_{h=1}^H{\log h\o h}-
{1\o 2}\log^2H\Bigr). \tag 6
\end{align*}

\vskip.3in
\noindent{\bf 2. Preliminaries}

\noindent{\bf Proposition 2.} For the given integer $a\geq 3$, let
$\chi$ be a Dirichlet character ${\rm mod}\,a$ and $\chi_0$ denote
the principal character. Then we have
\begin{align*}
\sum_{\substack{\chi({\rm mod}\,a)\\ \chi\ne\chi_0}}|L(1,\,\chi)|^2
&=\zeta(2)\varphi(a)\prod_{p|a}\Bigl(1-{1\o p^2}\Bigr)
-{\varphi^2(a)\o a^2}\Bigl(\log a+\sum_{p|a}{\log p\o p-1}\Bigr)^2\\
&\ +{\varphi^2(a)\o a^2}(\gamma^2+2\gamma_1-2\zeta(2))+
O\Bigl({\varphi(a)\o a^2}\Bigr),
\end{align*}
where $p$ denotes prime number, $\varphi(a)$ is the Euler totient
function, $\gamma$ is the Euler constant and $\gamma_1$ is defined
in (6).

This is Theorem 1 in [3].

\noindent{\bf Proposition 3.} For the given integer $a\geq 3$, we
have
$$
\sum_{\substack{\chi({\rm mod}\,a)\\ \chi(-1)=-1}}|L(1,\,\chi)|^2
={\pi^2\o 12}\cdot{\varphi^2(a)\o a^2}\Bigl(a\prod_{p|a}\Bigl(1+
{1\o p}\Bigr)-3\Bigr).
$$

One can see Theorem A in page 440 of [4].

\noindent{\bf Lemma 1.} In the neighborhood of $s=1$, there is
Laurent expansion
$$
\zeta^3(s)={1\o (s-1)^3}+{3\gamma\o (s-1)^2}+{3\gamma^2-3\gamma_1\o
s-1}+\cdots,
$$
where $\gamma$ is the Euler constant and $\gamma_1$ is defined in
(6).

\noindent{\bf Proof.} We know
$$
\zeta(s)={1\o s-1}+\gamma-\gamma_1(s-1)+\cdots.
$$
Hence,
\begin{align*}
\zeta^3(s)&=\Bigl({1\o s-1}+\gamma-\gamma_1(s-1)+\cdots\Bigr)^2
\Bigl({1\o s-1} +\gamma-\gamma_1(s-1)+\cdots\Bigr)\\
&=\Bigl({1\o (s-1)^2}+\gamma^2+{2\gamma\o s-1}-2\gamma_1+\cdots\Bigr)\\
&\ \cdot\Bigl({1\o s-1}+\gamma-\gamma_1(s-1)+\cdots\Bigr)\\
&=\Bigl({1\o (s-1)^2}+{2\gamma\o s-1}+(\gamma^2-2\gamma_1)+\cdots\Bigr)\\
&\ \cdot\Bigl({1\o s-1}+\gamma-\gamma_1(s-1)+\cdots\Bigr)\\
&={1\o (s-1)^3}+{\gamma\o (s-1)^2}-{\gamma_1\o s-1}+{2\gamma\o
(s-1)^2}+{2\gamma^2\o s-1}+{\gamma^2-2\gamma_1\o s-1}+\cdots\\
&={1\o (s-1)^3}+{3\gamma\o (s-1)^2}+{3\gamma^2-3\gamma_1\o
s-1}+\cdots.
\end{align*}

\noindent{\bf Lemma 2.} If $a\geq 3$, then
$$
\sum_{p|a}{\log p\o p-1}\ll\log\log a.
$$

\noindent{\bf Proof.} It is enough to prove for sufficiently large
$a$. When $x\geq\log a$, the function ${\log x\o x}$ decreases
monotonously. Thus
\begin{align*}
&\ \,\sum_{p|a}{\log p\o p-1}\ll\sum_{p|a}{\log p\o p}\\
&=\sum_{\substack{p |a\\ p\leq\log a}}{\log p\o p}+\sum_{\substack
{p|a\\ \log a< p}}{\log p\o p}\\
&\leq\sum_{p\leq\log a}{\log p\o p}+{\log\log a\o \log a}\sum_{p|a}1\\
&\ll\log\log a+{\log\log a\o \log a}\cdot\log a\\
&\ll\log\log a.
\end{align*}

\noindent{\bf Lemma 3.} If $a\geq 3$, then
$$
\sum_{p|a}{3p^3+2p^2+ 3p\o (p^2-1)^2}\cdot\log^2 p\ll\log\log^2a.
$$

\noindent{\bf Proof.} Assume that $a$ is sufficiently large. When
$x\geq\log a$, the function ${\log^2x\o x}$ decreases monotonously.
Thus
\begin{align*}
&\ \,\sum_{p|a}{3p^3+2p^2+ 3p\o (p^2-1)^2}\cdot\log^2 p\ll\sum_{p|a}
{\log^2p\o p}\\
&=\sum_{\substack{p |a\\ p\leq\log a}}{\log^2p\o p}+\sum_{\substack
{p|a\\ \log a< p}}{\log^2p\o p}\\
&\leq\log\log a\sum_{p\leq\log a}{\log p\o p}+{\log\log^2a\o \log a}
\sum_{p|a}1\\
&\ll\log\log^2a.
\end{align*}

\vskip.3in
\noindent{\bf 3. The proof of Theorem}

According to the discussion in [2], we have
\begin{align*}
S(N;\,a)={1\o \varphi(a)}\sum_{\chi({\rm mod}\,a)}\chi(-1){\rm
Res}_{s=1}\Bigl(f_\chi(s){N^s\o s}\Bigr)+\Delta(N;\,a), \tag 7
\end{align*}
where
\begin{align*}
f_\chi(s)={L(s,\,\chi_0)\o L(2s,\,\chi_0)}\cdot
L(s,\,\chi)L(s,\,\bar{\chi}) \tag 8
\end{align*}
and
$$
\Delta(N;\,a)\ll N^{1\o 2}\log^5N\cdot{a\o \varphi(a)}\prod_{p|a}
\Bigl(1-{1\o p^{1\o 2}}\Bigr)^{-1}.
$$
Write
\begin{align*}
S_1(N;\,a)={1\o \varphi(a)}\sum_{\substack{\chi({\rm mod}\,a)\\
\chi\ne\chi_0}}\chi(-1){\rm Res}_{s=1}\Bigl(f_\chi(s){N^s\o s}\Bigr)
\tag 9
\end{align*}
and
\begin{align*}
S_2(N;\,a)={1\o \varphi(a)}{\rm Res}_{s=1}\Bigl(f_{\chi_0}(s){N^s\o
s}\Bigr). \tag {10}
\end{align*}

The discussion in [2] shows that when $\chi\ne\chi_0$,
\begin{align*}
{\rm Res}_{s=1}\Bigl(f_\chi(s){N^s\o s}\Bigr)&={\rm
Res}_{s=1}\Bigl({L(s,\,\chi_0)L(s,\,\chi)L(s,\,\bar{\chi})N^s\o
L(2s,\,\chi_0)s}\Bigr)\\
&={6N\o \pi^2}\prod_{p|a}{p\o p+1}\cdot|L(1,\,\chi)|^2.
\end{align*}
Hence,
$$
S_1(N;\,a)={6N\o \pi^2}\cdot{1\o \varphi(a)}\prod_{p|a}{p\o p+1}
\cdot\sum_{\substack{\chi({\rm mod}\,a)\\ \chi\ne\chi_0}}\chi(-1)
|L(1,\,\chi)|^2.
$$

Since $\chi(-1)=1$ or $-1$, by Propositions 2 and 3, we get
\begin{align*}
&\ \,\sum_{\substack{\chi({\rm mod}\,a)\\ \chi\ne\chi_0}}\chi(-1)
|L(1,\,\chi)|^2\\
&=\sum_{\substack{\chi({\rm mod}\,a)\\
\chi\ne\chi_0\\ \chi(-1)=1}}|L(1,\,\chi)|^2-\sum_{\substack
{\chi({\rm mod}\,a)\\ \chi\ne\chi_0\\ \chi(-1)=-1}}|L(1,\,\chi)|^2\\
&=\sum_{\substack{\chi({\rm mod}\,a)\\ \chi\ne\chi_0}}
|L(1,\,\chi)|^2-2\sum_{\substack{\chi({\rm mod}\,a)\\ \chi\ne\chi_0
\\ \chi(-1)=-1}}|L(1,\,\chi)|^2\\
&=\sum_{\substack{\chi({\rm mod}\,a)\\ \chi\ne\chi_0}}
|L(1,\,\chi)|^2-2\sum_{\substack{\chi({\rm mod}\,a)\\ \chi(-1)=-1}}
|L(1,\,\chi)|^2\\
&=\zeta(2)\varphi(a)\prod_{p|a}\Bigl(1-{1\o p^2}\Bigr)-
{\varphi^2(a)\o a^2}\Bigl(\log a+\sum_{p|a}{\log p\o p-1}\Bigr)^2\\
&\ +{\varphi^2(a)\o a^2}(\gamma^2+2\gamma_1-2\zeta(2))+
O\Bigl({\varphi(a)\o a^2}\Bigr)\\
&\ -\zeta(2)\cdot{\varphi^2(a)\o a^2}\Bigl(a\prod_{p|a}\Bigl(1+{1\o
p}\Bigr)-3\Bigr)\\
&=\zeta(2)\varphi(a)\prod_{p|a}\Bigl(1-{1\o
p^2}\Bigr)-{\varphi^2(a)\o
a^2}\Bigl(\log a+\sum_{p|a}{\log p\o p-1}\Bigr)^2\\
&\ +{\varphi^2(a)\o a^2}(\gamma^2+2\gamma_1-2\zeta(2))+
O\Bigl({\varphi(a)\o a^2}\Bigr)\\
&\ -\zeta(2)\varphi(a)\prod_{p|a}\Bigl(1-{1\o p^2}\Bigr)+
3\zeta(2)\cdot{\varphi^2(a)\o a^2}\\
&=-{\varphi^2(a)\o a^2}\Bigl(\log a+\sum_{p|a}{\log p\o p-1}
\Bigr)^2+{\varphi^2(a)\o a^2}(\gamma^2+2\gamma_1+\zeta(2))
+O\Bigl({\varphi(a)\o a^2}\Bigr).
\end{align*}
Hence,
\begin{align*}
S_1(N;\,a)&={6N\o \pi^2}\cdot{\varphi(a)\o a^2}\prod_{p|a}{p\o p+1}
\cdot\Bigl(-\Bigl(\log a+\sum_{p|a}{\log p\o p-1}\Bigr)^2\\
&\ +(\gamma^2+2\gamma_1+\zeta(2))+O\Bigl({1\o \varphi(a)}
\Bigr)\Bigr)
\end{align*}
\begin{align*}
&={3\o \pi^2a}\prod_{p|a}{p-1\o p+1}\cdot N\Bigl(-2\log^2a-
4\log a\sum_{p|a}{\log p\o p-1}\\
&\ -2\Bigl(\sum_{p|a}{\log p\o p-1}\Bigr)^2+2\gamma^2+4\gamma_1
+2\zeta(2)+O\Bigl({1\o \varphi(a)}\Bigr)\Bigr).
\end{align*}

Now we proceed to calculate
$$
{\rm Res}_{s=1}\Bigl(f_{\chi_0}(s){N^s\o s}\Bigr).
$$
The discussion in [2] yields
\begin{align*}
f_{\chi_0}(s)={L^3(s,\,\chi_0)\o L(2s,\,\chi_0)}={\zeta^3(s)\o
\zeta(2s)}\prod_{p|a}{(p^s-1)^2\o p^s(p^s+1)}={\zeta^3(s)\o
\zeta(2s)}\cdot G(s), \tag{11}
\end{align*}
where
\begin{align*}
G(s)=\prod_{p|a}{(p^s-1)^2\o p^s(p^s+1)}. \tag {12}
\end{align*}

By Lemma 1,
\begin{align*}
&\ \,{\rm Res}_{s=1}\Bigl(f_{\chi_0}(s){N^s\o s}\Bigr)\\
&={\rm Res}_{s=1}\zeta^3(s)\cdot{N^s\o \zeta(2s)s}\cdot G(s)\\
&={\rm Res}_{s=1}{1\o (s-1)^3}\cdot{N^s\o \zeta(2s)s}\cdot G(s)\\
&\ + 3\gamma{\rm Res}_{s=1}{1\o (s-1)^2}\cdot{N^s\o
\zeta(2s)s}\cdot G(s)\\
&\ +(3\gamma^2-3\gamma_1){\rm Res}_{s=1}{1\o s-1}\cdot{N^s\o
\zeta(2s)s}\cdot G(s).
\end{align*}
We shall calculate these residues respectively.

1. It is easy to see
\begin{align*}
&\ \,(3\gamma^2-3\gamma_1){\rm Res}_{s=1}{1\o s-1}\cdot{N^s\o
\zeta(2s)s} \cdot G(s)\\
&=N\cdot{G(1)\o \zeta(2)}\cdot(3\gamma^2-3\gamma_1).
\end{align*}

2. We have
\begin{align*}
&\ \,3\gamma{\rm Res}_{s=1}{1\o (s-1)^2}\cdot{N^s\o \zeta(2s)s}
\cdot G(s)\\
&=3\gamma\Bigl({N^s\o \zeta(2s)s}\cdot G(s)\Bigr)'\Bigl|_{s=1}\\
&=3\gamma\Bigl(\Bigl({N^s\o \zeta(2s)s}\Bigr)'\cdot G(s)+{N^s\o
\zeta(2s)s}\cdot G'(s)\Bigr)\Bigl|_{s=1}\\
&=3\gamma\Bigl({N^s\log N\cdot\zeta(2s)s-N^s(\zeta'(2s)\cdot 2s
+\zeta(2s))\o (\zeta(2s)s)^2}\cdot G(s)\\
&\ +{N^s\o \zeta(2s)s}\cdot G'(s)\Bigr)\Bigl|_{s=1}\\
&=3\gamma G(1)\Bigl({N\log N\cdot\zeta(2)-N(2\zeta'(2)+\zeta(2))\o
\zeta^2(2)}+{N\o \zeta(2)}\cdot {G'(1)\o G(1)}\Bigr)\\
&=N\cdot{G(1)\o \zeta(2)}\cdot3\gamma\Bigl({\zeta(2)\log N
-(2\zeta'(2)+\zeta(2))\o \zeta(2)}+{G'(1)\o G(1)}\Bigr)\\
&=N\cdot{G(1)\o \zeta(2)}\cdot\Bigl(3\gamma\log
N-3\gamma\Bigl(2{\zeta'(2)\o \zeta(2)}+1\Bigr)+3\gamma\cdot{G'(1)\o
G(1)}\Bigr).
\end{align*}

The derivative of $\log G(s)$ is
\begin{align*}
{G'(s)\o G(s)}&=\sum_{p|a}(2\log(p^s-1)-s\log p-\log(p^s+1))'\\
&=\sum_{p|a}\Bigl({2\o p^s-1}\cdot p^s\log p-\log p-{1\o p^s+1}\cdot
p^s\log p\Bigr) \tag {13}\\
&=\sum_{p|a}\Bigl({2\o p^s-1}+{1\o p^s+1}\Bigr)\log p.
\end{align*}
Hence,
\begin{align*}
{G'(1)\o G(1)}=\sum_{p|a}\Bigl({2\o p-1}+{1\o p+1}\Bigr)\log
p=\sum_{p|a}{3p+1\o p^2-1}\cdot\log p, \tag {14}
\end{align*}
which is a formula in page 1652 of [2].

3. We have
\begin{align*}
&\ \,{\rm Res}_{s=1}{1\o (s-1)^3}\cdot{N^s\o \zeta(2s)s}\cdot G(s)\\
&={1\o 2!}\Bigl({N^s\o \zeta(2s)s}\cdot G(s)\Bigr)''\Bigl|_{s=1}
\end{align*}
\begin{align*}
&={1\o 2}\Bigl(\Bigl({N^s\o \zeta(2s)s}\Bigr)''\cdot
G(s)+2\Bigl({N^s\o \zeta(2s)s}\Bigr)'\cdot G'(s)+{N^s\o \zeta(2s)s}
\cdot G''(s)\Bigr)\Bigl|_{s=1}\\
&={1\o 2}\cdot\Bigl({N^s\o \zeta(2s)s}\Bigr)''\cdot G(s)\Bigl|_{s=1}
+\Bigl({N^s\o \zeta(2s)s}\Bigr)'\cdot G'(s)\Bigl|_{s=1}\\
&\ +{1\o 2}\cdot{N^s\o \zeta(2s)s}\cdot G''(s) \Bigl|_{s=1}.
\end{align*}
We shall calculate these expressions respectively.

a) Since
\begin{align*}
&\ \,\Bigl({N^s\o \zeta(2s)s}\Bigr)''=\Bigl({N^s\log
N\cdot\zeta(2s)s-N^s(\zeta'(2s)\cdot 2s +\zeta(2s))\o (\zeta(2s)
s)^2}\Bigr)'\\
&={(N^s\log N\cdot\zeta(2s)s-N^s(\zeta'(2s)\cdot 2s+\zeta(2s)))'
(\zeta(2s)s)^2\o (\zeta(2s)s)^4}\\
&\ -{(N^s\log N\cdot\zeta(2s)s-N^s(\zeta'(2s)\cdot 2s+\zeta(2s)))
2\zeta(2s)s(\zeta(2s)s)'\o (\zeta(2s)s)^4}\\
&={(N^s\log N\cdot\zeta(2s)s)'\o (\zeta(2s)s)^2}-{(N^s
(\zeta'(2s)\cdot 2s+\zeta(2s)))'\o (\zeta(2s)s)^2}\\
&\ -{N^s\log N\cdot\zeta(2s)s\cdot 2(\zeta'(2s)\cdot 2s+
\zeta(2s))\o (\zeta(2s)s)^3}\\
&\ +{N^s(\zeta'(2s)\cdot 2s+\zeta(2s))\cdot 2(\zeta'(2s)\cdot
2s+\zeta(2s))\o (\zeta(2s)s)^3}\\
&={(N^s\log^2N\cdot\zeta(2s)s+N^s\log N\cdot\zeta'(2s)2s+N^s\log
N\cdot\zeta(2s))\o (\zeta(2s)s)^2}\\
&\ -{N^s\log N(\zeta'(2s)\cdot 2s+\zeta(2s))+N^s(\zeta''(2s)\cdot
4s+2\zeta'(2s)+2\zeta'(2s))\o (\zeta(2s)s)^2}\\
&\ -{N^s\log N\cdot 2(\zeta'(2s)\cdot 2s+\zeta(2s))\o
(\zeta(2s)s)^2}+{N^s\cdot 2(2s\zeta'(2s)+\zeta(2s))^2\o
(\zeta(2s)s)^3},
\end{align*}
\begin{align*}
&\ \,{1\o 2}\cdot\Bigl({N^s\o \zeta(2s)s}\Bigr)''\cdot G(s)\Bigl|_{s=1}\\
&=N\cdot{G(1)\o 2}\cdot\Bigl({\zeta(2)\log^2N+(2\zeta'(2)
+\zeta(2))\log N\o \zeta^2(2)}\\
&\ -{(2\zeta'(2)+\zeta(2))\log N+4\zeta''(2)+4\zeta'(2)\o
\zeta^2(2)}\\
&\ -{2(2\zeta'(2)+\zeta(2))\log N\o \zeta^2(2)}+{2(2\zeta'(2)
+\zeta(2))^2\o \zeta^3(2)}\Bigr)\\
&=N\cdot{G(1)\o \zeta(2)}\cdot\Bigl({1\o 2}\log^2N-{1\o
\zeta(2)}(2\zeta'(2)+\zeta(2))\log N\\
&\ +{1\o \zeta^2(2)}(2\zeta'(2)+\zeta(2))^2-{2\o \zeta(2)}
(\zeta''(2)+\zeta'(2))\Bigr).
\end{align*}

b) We have
\begin{align*}
&\ \,\Bigl({N^s\o \zeta(2s)s}\Bigr)'\cdot G'(s)\Bigl|_{s=1}\\
&={N^s\log N\cdot\zeta(2s)s-N^s (\zeta'(2s)\cdot 2s +\zeta(2s))\o
(\zeta(2s)s)^2}\cdot G'(s)\Bigl|_{s=1}\\
&={N\log N\cdot\zeta(2)-N(2\zeta'(2)+\zeta(2))\o \zeta^2(2)}\cdot
G'(1)\\
&=N\cdot{G(1)\o \zeta(2)}\cdot{\zeta(2)\log N-(2\zeta'(2)+
\zeta(2))\o \zeta(2)}\cdot {G'(1)\o G(1)}\\
&=N\cdot{G(1)\o \zeta(2)}\cdot\Bigl(\log N-2{\zeta'(2)\o
\zeta(2)}-1\Bigr)\cdot {G'(1)\o G(1)}.
\end{align*}

c) Differentiating the equality in (13)
$$
G'(s)=G(s)\sum_{p|a}\Bigl({2\o p^s-1}+{1\o p^s+1}\Bigr)\log p,
$$
we get
\begin{align*}
G''(s)&=G'(s)\sum_{p|a}\Bigl({2\o p^s-1}+{1\o p^s+1}\Bigr)\log p\\
&\ +G(s)\sum_{p|a}\Bigl(-{2p^s\log p\o (p^s-1)^2}-{p^s\log p\o
(p^s+1)^2}\Bigr)\log p\\
&=G(s)\Bigl({G'(s)\o G(s)}\sum_{p|a}\Bigl({2\o p^s-1}+{1\o p^s+1}
\Bigr)\log p\\
&\ -\sum_{p|a}\Bigl({2\o (p^s-1)^2}+{1\o (p^s+1)^2}\Bigr)
p^s\log^2p\Bigr).
\end{align*}
Therefore
\begin{align*}
G''(1)&=G(1)\Bigl({G'(1)\o G(1)}\sum_{p|a}\Bigl({2\o p-1}+{1\o
p+1}\Bigr)\log p\\
&\ -\sum_{p|a}\Bigl({2\o (p-1)^2}+{1\o (p+1)^2}\Bigr)p\log^2p\Bigr)\\
&=G(1)\Bigl(\Bigl({G'(1)\o G(1)}\Bigr)^2-\sum_{p|a}{3p^3+2p^2+3p\o
(p^2-1)^2}\cdot\log^2p\Bigr),
\end{align*}
which is a formula in page 1652 of [2].

We have
\begin{align*}
&\ \,{1\o 2}\cdot{N^s\o \zeta(2s)s}\cdot G''(s)\Bigl|_{s=1}\\
&={N\o 2\zeta(2)}\cdot G''(1)\\
&=N\cdot{G(1)\o \zeta(2)}\cdot\Bigl({1\o 2}\Bigl({G'(1)\o
G(1)}\Bigr)^2 -{1\o 2}\sum_{p|a}{3p^3+2p^2+3p\o (p^2-1)^2}\cdot
\log^2p\Bigr).
\end{align*}

Combining results in cases a), b) and c), we get
\begin{align*}
&\ \,{\rm Res}_{s=1}{1\o (s-1)^3}\cdot{N^s\o \zeta(2s)s}\cdot G(s)\\
&=N\cdot{G(1)\o \zeta(2)}\cdot\Bigl(\Bigl\{{1\o 2}\log^2N-{1\o
\zeta(2)}(2\zeta'(2)+\zeta(2))\log N\\
&\ +{1\o \zeta^2(2)}(2\zeta'(2)+\zeta(2))^2-{2\o \zeta(2)}
(\zeta''(2)+\zeta'(2))\Bigr\}\\
&\ +\Bigl\{{G'(1)\o G(1)}\cdot\log N-\Bigl(2{\zeta'(2)\o
\zeta(2)}+1\Bigr)\cdot{G'(1)\o G(1)}\Bigr\}\\
&\ +\Bigl\{{1\o 2}\Bigl({G'(1)\o G(1)}\Bigr)^2-{1\o 2}
\sum_{p|a}{3p^3+2p^2+3p\o (p^2-1)^2}\cdot\log^2p\Bigr\}\Bigr)\\
&=N\cdot{G(1)\o \zeta(2)}\cdot\Bigl({1\o 2}\log^2N-{1\o \zeta(2)}
(2\zeta'(2)+\zeta(2))\log N\\
&\ +{G'(1)\o G(1)}\cdot\log N+{1\o 2}\Bigl({G'(1)\o G(1)}\Bigr)^2\\
&\ -{1\o 2}\sum_{p|a}{3p^3+2p^2+3p\o (p^2-1)^2}\cdot
\log^2p-\Bigl(2{\zeta'(2)\o \zeta(2)}+1\Bigr)\cdot{G'(1)\o G(1)}\\
&\ +{1\o \zeta^2(2)}(2\zeta'(2)+\zeta(2))^2-{2\o \zeta(2)}
(\zeta''(2)+\zeta'(2))\Bigr).
\end{align*}

Now we can get from cases 1, 2 and 3 that
\begin{align*}
&\ \,{\rm Res}_{s=1}\Bigl(f_{\chi_0}(s){N^s\o s}\Bigr)\\
&=N\cdot{G(1)\o \zeta(2)}\cdot\Bigl({1\o 2}\log^2N-{1\o \zeta(2)}
(2\zeta'(2)+\zeta(2))\log N\\
&\ +{G'(1)\o G(1)}\cdot\log N+{1\o 2}\Bigl({G'(1)\o G(1)}\Bigr)^2
\end{align*}
\begin{align*}
&\ -{1\o 2}\sum_{p|a}{3p^3+2p^2+3p\o (p^2-1)^2}\cdot
\log^2p-\Bigl(2{\zeta'(2)\o \zeta(2)}+1\Bigr)\cdot{G'(1)\o G(1)}\\
&\ +{1\o \zeta^2(2)} (2\zeta'(2)+\zeta(2))^2-{2\o \zeta(2)}
(\zeta''(2)+\zeta'(2))\Bigr)\\
&\ +N\cdot{G(1)\o \zeta(2)}\cdot\Bigl(3\gamma\log
N-3\gamma\Bigl(2{\zeta'(2)\o \zeta(2)}+1\Bigr)+3\gamma\cdot{G'(1)
\o G(1)}\Bigr)\\
&\ +N\cdot{G(1)\o \zeta(2)}\cdot(3\gamma^2-3\gamma_1)\\
&=N\cdot{G(1)\o \zeta(2)}\cdot\Bigl({1\o 2}\log^2N+3\gamma\log N\\
&\ -{1\o \zeta(2)}(2\zeta'(2)+\zeta(2))\log N+{G'(1)\o G(1)}\cdot
\log N\\
&\ +{1\o 2}\Bigl({G'(1)\o G(1)}\Bigr)^2-{1\o 2} \sum_{p|a}{3p^3
+2p^2+3p\o (p^2-1)^2}\cdot\log^2p\\
&\ +\Bigl(3\gamma-2{\zeta'(2)\o \zeta(2)}-1\Bigr)\cdot{G'(1)\o
G(1)}+{1\o \zeta^2(2)}(2\zeta'(2)+\zeta(2))^2\\
&\ -{2\o \zeta(2)}(\zeta''(2)+\zeta'(2))-3\gamma\Bigl(2{\zeta'(2)\o
\zeta(2)}+1\Bigr)+3\gamma^2-3\gamma_1\Bigr).
\end{align*}

Since
\begin{align*}
{1\o \varphi(a)}\cdot{G(1)\o \zeta(2)}&={1\o a}\prod_{p|a}{p\o p-1}
\cdot{6\o \pi^2}\prod_{p|a}{(p-1)^2\o p(p+1)}\\
&={6\o \pi^2 a}\prod_{p|a}{p-1\o p+1},
\end{align*}
\begin{align*}
S_2(N;\,a)&={1\o \varphi(a)}{\rm Res}_{s=1}\Bigl(f_{\chi_0}(s)
{N^s\o s}\Bigr)\\
&={3\o \pi^2a}\prod_{p|a}{p-1\o p+1}\cdot N\Bigl(\log^2N+
\Bigl(6\gamma-4{\zeta'(2)\o \zeta(2)}-2\Bigr)\log N\\
&\ +\Bigl(\sum_{p|a}{6p+2\o p^2-1}\cdot\log p\Bigr)\log N
+\Bigl(\sum_{p|a}{3p+1\o p^2-1}\cdot\log p\Bigr)^2\\
&\ -\sum_{p|a}{3p^3+2p^2+3p\o (p^2-1)^2}\cdot\log^2p\\
&\ +\Bigl(6\gamma-4{\zeta'(2)\o \zeta(2)}-2\Bigr)\sum_{p|a}{3p+1\o
p^2-1}\cdot\log p
\end{align*}
\begin{align*}
&\ +{2\o \zeta^2(2)}(2\zeta'(2)+\zeta(2))^2-{4\o \zeta(2)}
(\zeta''(2)+\zeta'(2))\\
&\ -6\gamma\Bigl(2{\zeta'(2)\o \zeta(2)}+1\Bigr)+6\gamma^2
-6\gamma_1\Bigr).
\end{align*}

Therefore
\begin{align*}
&\ \,S_1(N;\,a)+S_2(N;\,a)\\
&={3\o \pi^2a}\prod_{p|a}{p-1\o p+1}\cdot N\Bigl(-2\log^2a-4\log
a\sum_{p|a} {\log p\o p-1}\\
&\ -2\Bigl(\sum_{p|a}{\log p\o p-1}\Bigr)^2+2\gamma^2+4\gamma_1
+2\zeta(2)+O\Bigl({1\o \varphi(a)}\Bigr)\Bigr)\\
&\ +{3\o \pi^2a}\prod_{p|a}{p-1\o p+1}\cdot N\Bigl(\log^2N+
\Bigl(6\gamma-4{\zeta'(2)\o \zeta(2)}-2\Bigr)\log N\\
&\ +\Bigl(\sum_{p|a}{6p+2\o p^2-1}\cdot\log p\Bigr)\log N+
\Bigl(\sum_{p|a}{3p+1\o p^2-1}\cdot\log p\Bigr)^2\\
&\ -\sum_{p|a}{3p^3+2p^2+3p\o (p^2-1)^2}\cdot\log^2p\\
&\ +\Bigl(6\gamma-4{\zeta'(2)\o \zeta(2)}-2\Bigr)\sum_{p|a}
{3p+1\o p^2-1}\cdot\log p\\
&\ +{2\o \zeta^2(2)}(2\zeta'(2)+\zeta(2))^2-{4\o \zeta(2)}
(\zeta''(2)+\zeta'(2))\\
&\ -6\gamma\Bigl(2{\zeta'(2)\o\zeta(2)}+1\Bigr)+6\gamma^2-
6\gamma_1\Bigr)\\
&={3\o \pi^2a}\prod_{p|a}{p-1\o p+1}\cdot N\Bigl(\log^2N
+\Bigl(6\gamma-4{\zeta'(2)\o \zeta(2)}-2\Bigr)\log N\\
&\ +\Bigl(\sum_{p|a}{6p+2\o p^2-1}\cdot\log p\Bigr)\log N-
2\log^2a-4\log a\sum_{p|a} {\log p\o p-1}\\
&\ -2\Bigl(\sum_{p|a}{\log p\o p-1}\Bigr)^2+\Bigl(\sum_{p|a}{3p+1\o
p^2-1}\cdot\log p\Bigr)^2\\
&\ -\sum_{p|a}{3p^3+2p^2+3p\o (p^2-1)^2}\cdot\log^2p\\
&\ +\Bigl(6\gamma-4{\zeta'(2)\o \zeta(2)}-2\Bigr)\sum_{p|a}{3p+1\o
p^2-1}\cdot\log p\\
&\ +{2\o \zeta^2(2)}(2\zeta'(2)+\zeta(2))^2-{4\o \zeta(2)}
(\zeta''(2)+\zeta'(2))
\end{align*}
\begin{align*}
&\ -6\gamma\Bigl(2{\zeta'(2)\o\zeta(2)}+1\Bigr)+8\gamma^2-
2\gamma_1+2\zeta(2)+O\Bigl({1\o \varphi(a)}\Bigr)\Bigr)\\
&={3\o \pi^2a}\prod_{p|a}{p-1\o p+1}\cdot N(\log^2N+c_1(a) \log
N+c_0(a)),
\end{align*}
where $c_1(a)$ and $c_0(a)$ are defined in (2) and (5). By Lemmas 2
and 3, we can see that terms of $c_0(a)$ in (5) are arranged in
decreasing order of $a$.

So far the proof of Theorem is finished.

\vskip.3in
\noindent{\bf Acknowledgements.}\ \ This work was
supported by National Natural Science Foundation of China (Grant No.
11071235). I would like to thank Professor Hao Pan of Nanjing
University for telling me the publication of paper [2].


\bigskip

\

Institute of  Mathematics, Academia Sinica, Beijing 100190, P. R.
China

E-mail: jiach@math.ac.cn

\vskip.6in

\begin{thebibliography}{9}

\bibitem{1} R. K. Guy, {\it Unsolved Problems in Number Theory},
third edition, Springer-Verlag, 2004.

\bibitem{2} Jingjing Huang and R. C. Vaughan, {\it Mean value
theorems for binary Egyptian fractions}, J. Number Theory, {\bf 131}
(2011), 1641-1656.

\bibitem{3} S. Kanemitsu, Y. Tanigawa, M. Yoshimoto and Wenpeng Zhang,
{\it On the discrete mean square of Dirichlet $L-$functions at 1},
Math. Z., {\bf 248}(2004), 21-44.

\bibitem{4} Wenpeng Zhang, {\it On the mean values of Dedekind sums},
Journal de Th\'eorie des Nombres de Bordeaux, {\bf 8}(1996),
429-442.






\end{thebibliography}
\end{document}